# Chinese Names for Integers

**Rémi ANICOTTE**

CRLAO, France

Abstract: Chinese names for integers have always used the digits [1] through [9] and a series of decimal pivots starting with [10], [$10^2$], [$10^3$] and [$10^4$]. Changes occurred in the way the compounds [*digit*][*pivot*] were concatenated, with the conjunction *yòu* until the 3rd century BCE, then with the term *líng*, which emerged around the 12th century CE. The behavior of the morpheme [1] with pivots also evolved. Finally, in Contemporary Chinese, there is a choice between two morphemes for the digit 2 yielding legitimate alternative numerals; and there is the possibility to form elliptic number names which are not meant to be incorporated before classifiers. Some changes in the features of Chinese linguistic numeration were likely the result of language planning; they nevertheless hint at a tension between a tendency to maintain the morphosyntax of number names within the framework of the syntax of quantification versus simplification and shorter numerals.

List of Abbreviations: CLF: classifier; MW: measure word; Num: numeral; PART: particle; PL: plural; 3OBJ; 3SG: third person singular pronoun; [*n*] (with a number *n* written in Arabic digits): the mono-morphemic expression of the number *n* in a given language; A(B): the character A is a rendition of the original character encountered in the Chinese corpus, the character B in parenthesis is a modern form for what A is understood to mean.

Acknowledgements: For their invaluable help in selecting relevant corpus and/or interpreting available data, I thank Karine Chemla (SPHERE), Redouane Djamouri (CRLAO), Christoph Harbsmeier (University of Oslo), Laurent Sagart (CRLAO), Sylviane Schwer (Paris 13 University), Xu Dan (INaLCO, IUF) and Zhang Xiancheng (Southwest University, Chongqing). Conclusions and shortcomings are mine.

# 1. DIGITS AND PIVOTS IN CHINESE NUMBER NAMES

Miller et al. (1995, 2005) commented that Chinese names[1] for numbers between 11 and 99 have a closer association with the positional notation in Arabic digits for Chinese than their English equivalents because there are no special words for teens and tens[2] in Chinese[3]. Miller and al.'s comparison was limited to numbers from 1 to 99 in Chinese and English. But for numbers over 100, Contemporary Chinese has number names which do not map well onto the positional notation of numbers, and even has free variants for some numbers in some situations.

A straightforward method to represent number names is to write down a linear string of signs, each symbolizing one morpheme in the order of speech production. I use the notation [*number*] with a number written in Arabic digits between square brackets to represent the mono-morphemic item which expresses the bracketed number in a given language. For example the notation [10] represents *ten* in English as well as *shí* in Chinese. The notation [$10^4$] represents *wàn* in Chinese, but would not occur in the representation of English number names because 10,000 is expressed as *ten thousand,* which we symbolize as [10][$10^3$], a compound of the mono-morphemic items *ten* [10] and *thousand* [$10^3$]. This representation accounts for the linearity of speech and notes all and only what is said[4]. Whereas arithmetical translations such as "1 x $10^2$ + 3 x 10 + 1" for [1][$10^2$][3][10][1] can be handy on occasions, they fail to give a proper account of the linguistic sequence, as there are usually no linguistic counterparts to the symbols for multiplication and addition, and generally no surface differences between the various numerical morphemes.

Number names are sequences of numerical morphemes and occasionally of linking words. Among the numerical morphemes, we need to differentiate between *digits* and *pivots,* which manifest semantic and syntactic disparities.

I use the name *multiplicative pivots* for numerical morphemes (and also the corresponding numbers) which are used to build the names of at least some of the numbers that are multiples of the pivot's value. For example, *hundred* [$10^2$] is a multiplicative pivot in English found in *one hundred* [10$^2$], *two hundred* [2][10$^2$], *three hundred* [3][10$^2$], etc [5]. The words *two* [2] and *thirty* [30] [6] are not multiplicative pivots because they are not involved in forming the English names of their multiples 4, 6, etc. and 60, 90, etc.

---

[1] I use the words *number names* and *numerals* for expressions produced by *linguistic numerations*.

[2] Miller and al. (2005) addressed the better performance of Chinese children compared to US children in numeracy and claimed that the characteristics of Chinese number names is an advantage for the preschool learning of numeracy. The authors explained that the performance gap was also related to cultural differences in attitudes towards education.

[3] A discussion of the mono-morphemic expressions for 20 in modern Hakka and Cantonese, or arguably for 20 and other tens in the history of standard Chinese is beyond the scope of this paper.

[4] Brainerd (1968) and Brainerd & Peng (1968) represented these strings of morphemes with a succession of numbers in Arabic digits separated by blanks, for example the latter article represents *yī bǎi sān shí yī* (i.e. the name for 131 in Contemporary Chinese) as # 1 100 3 10 1 # while I write [1][$10^2$][3][10][1].

Tab. 1: Linguistic decimal scales of American English and Contemporary Chinese

| Decimal ranks | American English | Contemporary Chinese |
|---|---|---|
| 10 | *ten* (not a pivot) | *shí* 十 (a pivot) |
| 10$^2$ | *hundred* | *bǎi* 百 |
| 10$^3$ | *thousand* | *qiān* 千 |
| 10$^4$ | | *wàn* 萬 |
| 10$^5$ | | |
| 10$^6$ | *million* | |
| 10$^7$ | | |
| 10$^8$ | | *yì* 億 |
| 10$^9$ | *billion* | |

Contemporary Chinese multiplicative pivots are *shí* 十 [10], *bǎi* 百 [10$^2$], *qiān* 千 [10$^3$], *wàn* 萬 [10$^4$] and *yì* 億 [10$^8$][7]; they are monomorphemic names of units taken in a decimal scale. An *arithmetical* decimal scale is made of the series of powers of 10, but a *linguistic* decimal scale in a given language does not need to cover the whole arithmetical scale. In American English the common decimal scale pivots are *hundred* [$10^2$], *thousand* [$10^3$], *million* [$10^6$] and *billion* [$10^9$][8], with only $10^2$ and $10^3$ shared with Chinese; this is shown in Table 1 and proves that the same strategy of resorting to pivots on the same arithmetical scale does not necessarily imply an identical choice of which ranks have a monomorphemic name. The highest pivot before a gap in the series is [$10^3$] in English and [$10^4$] in Chinese. They are the first *outer* pivots, the previous ones being *inner* pivots[9].

---

[5] This definition allows for example French *cent* [10$^2$] to be classified as a multiplicative pivot because *cent* [10$^2$] is found in *deux cents* [2][10$^2$]PL, *trois cents* [3][10$^2$]PL, etc. and even though the name for 100 is simply *cent* [10$^2$] which does not call for the digit *un* [1].

[6] English *thirty* [30] is an *additive pivot* used only additively to form the number names *thirty one* [30][1] up to *thirty nine* [30][9].

[7] These are the sole multiplicative pivots taught in today's primary and secondary school curricula. The history of the evolution of the list of pivots beyond *wàn* [10$^4$] is beyond the scope of this paper. It ended with the standardization of a pivot *yì* [10$^8$] and no agreed forms above. The word *zhào* 兆 is often cited, but no one value was ever agreed on: there are three contradictory definitions. The value 10$^6$ is the equivalent of the prefix *mega* of the international system of units, and is the one accepted in today's curricula and the only one mentioned in the 2010 edition of the dictionary *Xiàndài hànyǔ guīfàn cídiǎn* 现代汉语规范词典 published in Beijing. The values 10$^{12}$ or 10$^{16}$ are also attested; they are *wàn yì* [10$^4$][10$^8$] and *yì yì* [10$^8$][10$^8$] respectively and make *zhào* a regular (but contradicting) extension of the series of pivots.

I use the name *digits* for numerical morphemes (and the corresponding numbers) which are used additively with pivots to form the names of consecutive numbers or multiplicatively to form the names of consecutive multiples. The Chinese digits are *yī* 一 [1][10], *èr* 二 [2], *sān* 三 [3], *sì* 四 [4], *wǔ* 五 [5], *liù* 六 [6], *qī* 七 [7], *bā* 八 [8], *jiǔ* 九 [9][11].

## 2. HISTORICAL SURVEY OF CHINESE NUMBER NAMES

The data available on Chinese linguistic numeration is spread over a period of three thousand years. The number name system was always decimal and relied on the same digits and scale of multiplicative pivots starting with [10], [$10^2$], [$10^3$] and [$10^4$]. But changes occurred in the manner the compounds [*digit*][*pivot*] were juxtaposed or concatenated using the terms *yòu* and *líng*. The way the morpheme [1] was used with pivots also evolved. Finally a choice between two morphemes to express the digit 2 in the names of exact numbers emerged in Contemporary Chinese.

---

[8] The first decimal rank name in English is *ten* [10] which is not a multiplicative pivot because tens are not compounds of the morpheme *ten* [10], they are expressed, putting aside the etymological level of analysis, with the mono-morphemic words *twenty* [20], *thirty* [30], etc.

[9] The importance of this difference was pointed out by Sylviane Schwer (Paris 13 University). This is not a terminological quibble since some features of Chinese number names depend on it.

[10] According to the rules of *pīnyīn* transcription, the digit [1] is always Romanized *yī* with a first tone mark regardless of the actual tone in Contemporary Chinese. This tone depends on that of the following syllable; for example [1][$10^4$] is pronounced *yí wàn* what is of interest only in publications concerned with pronunciation.

[11] English digits are also 1 to 9: the corresponding morphemes are used additively for example in *twenty one* [20][1] to *twenty nine* [20][9], and multiplicatively to form for example *one hundred* [1][$10^2$] to *nine hundred* [9][$10^2$]. But English *ten* [10] to *nineteen* [19], although all mono-morphemic, if one accepts disregarding the etymological level of analysis, are neither digits nor pivots.

## *2-1. The linking terms* **yòu** *and* **líng**

Most ancient excavated data dates back to Shang inscriptions on oracle bones and bronze vessels (13th to 11th centuries BCE). Then there are Zhou (11th to 5th centuries BCE) and Warring States (5th to 3rd century BCE) inscriptions on bronze vessels[12].

An important feature which can be observed in Shang and Zhou inscriptions is that the conjunction *yòu* 有(又)[13] was sometimes used to link tens and units, and sometimes also hundreds and tens. But this was not obligatory; for example both *shí yòu wǔ*, i.e. [10] *yòu* [5], and *shí wǔ*, i.e. [10][5], are encountered. The use of a conjunction like *yòu* keeps the morphosyntax of number names rooted in the syntax of noun phrases, whereas a linguistic numeration can well be consistent without linking terms interrupting the chain of compounds [*digit*][*pivot*].

In Shang inscriptions on oracle bones, expressions using *yòu* between tens and units represented only 5% of all expressions involving tens and units (59 instances with *yòu* versus 1175 without it). In the available bronze inscriptions from the Zhou dynasty, the proportion reaches approximately 98% (there are 284 instances with *yòu* versus 5 without it). This discrepancy could reflect a genuine linguistic evolution or a mere stylistic difference: oracular inscriptions could be more stenographic than pompous inscriptions on bronze vessels, which would reflect the spoken pattern of officials.

In any case, later bronze vessel inscriptions dating to the Warring States period (5th to 3rd century BCE) reveal an indisputable linguistic change, because only around 8% of expressions still used *yòu* (24 instances with it versus 267 without it). This trend to discontinue use of *yòu* had already reached its full development at the beginning of the 2nd century BCE since the conjunction is no longer to be found in the names for integers written in the *Suàn shù shū*[14]. In this corpus, [*digit*][*pivot*] compounds are directly concatenated. The conjunction is used only in expressions for mixed numbers, in order to link an integer and a fraction smaller than one, which can be seen in (1).

---

[12] The data was accessed in CHANT on several occasions between December 2011 and November 2012.

[13] This notation 有(又) is conventional in publications on excavated Chinese texts. The character 有 is a rendition of the character encountered in the corpus, and the character 又 in parenthesis is the modern form of what is understood for the original character.

[14] The *Suàn shù shū* 算數書 was excavated from a tomb where a calendar for the year 186 BCE was found; so that the tomb is thought to have been closed that very year and the manuscript was probably written in the beginning of the 2nd century BCE. Peng Hao (2001: 4-6) states that some parts were copied from texts originally written in the kingdom of Qin before the unification of China in 221 BCE, while other texts could only have been composed during the reign of the Western Han dynasty which began in 206 BCE.

(1)       十六尺         有(又)
in *Suàn shù shū*   *shí liù chǐ*   *yòu*
strip 55       [10][6] *chǐ*    and
          '16 12/18 *chǐ*' (expressed as "16 *chǐ* and 12/18 *chǐ*";

         十八分       尺    之    十二
         *shí bā fēn*   *chǐ*   *zhī*   *shí èr*
         [10][8] part   *chǐ*   *zhī*   [10][2]
         *chǐ* is a unit of length)

Now we must turn to transmitted texts in order to investigate the use of the morpheme *líng* 零 in the Chinese number name system. According to Xu Pinfang and Zhang Hong (2006: 101), its first known appearance is in the word *èr bǎi líng qī*, i.e. [2][$10^2$] *líng* [7], for the number 207, found in a calendar[15] published in 1180 CE. Actually there is also a similar instance of the morpheme *dān* 單 found in the word *liù bǎi dān yī*, i.e. [6][$10^2$] *dān* [1], for the number 601 in writings[16] published in 1270 CE. The word *líng* originally meant *raindrops remaining on objects after a rainfall*; *dān* meant *alone, isolated*. They both introduce the remaining odd units of a number with hundreds but no tens; this manner of speech could have first emerged as a free construction. The word *lìng* 另 (*another*) is also encountered; it could be a mere graphical variation for *líng* 零. The term *líng* was the most common[17] and is the only one remaining in the number names of Contemporary Chinese. The numerical expressions formed with *líng* are no longer free variants but standard number names; however, the details of the process of standardization are unclear.

Qin Jiushao, the Song dynasty author of the *Shù shū jiǔ zhāng* (1247)[18], regularly used *líng* after [$10^4$], [$10^3$] and [$10^2$] when the following digit was not associated respectively with [$10^3$], [$10^2$] and [10]. Shi Yuechun did the same in his *Bǎi jī shù yǎn* (1861)[19]. But Li Zhizao in his *Tóng wén suàn zhǐ* (1613)[20], 366 years after Qin

---

[15] The *Dà míng lì* 大明歷 published by Zhao Zhiwei 趙知微.

[16] The *Zhū zǐ yǔ lèi* 朱子語類.

[17] In CCL (accessed in January 2012), there are 336 occurrences of [$10^2$] *líng* versus only 60 occurrences of [$10^2$] *dān* and 4 of [$10^2$] *lìng*; 69 occurrences of [$10^3$] *líng* versus 2 of [$10^3$] *dān* and 2 of [$10^3$] *lìng*; 45 occurrences of [$10^4$] *líng* versus 7 of [$10^4$] *dān* and 2 of [$10^4$] *lìng*.

[18] The book of mathematics *Shù shū jiǔ zhāng* 數書九章 was written by Qin Jiushao 秦九韶 (1202-1261), I checked an edition prepared in 1842.

[19] The *Bǎi jī shù yǎn* 百鸡术衍 was written by Shi Yuechun 时日醇 (1807-1880); I checked an edition from 1872.

[20] The *Tóng wén suàn zhǐ* 同文算指 was published in 1613 under the Ming dynasty. It was written by Li Zhizao 李之藻 (1565-1630), with plausible participation by Xu Guangqi 徐光啟 (1562-1633) and Matteo Ricci (1552-1610). I checked the *Sì kù quán shū* 四庫全書 edition prepared in the 18th century.

Jiushao and 248 years before Shi Yuechun, used *líng* only on occasions within integer names; for example it is used in the expression for 4004 in (2) but not for 2002.

(2)  四千零四分　　鏊　　之　二千一百三十
in *Tóng wén*　*sì qiān líng sì fēn*　*lí*　*zhī*　*èr qiān yī bǎi sān shí*
*suàn zhǐ*　[4][$10^3$] *líng* [4] *fēn*　hundredth　*zhī*　[2][$10^3$][1][$10^2$][3][10]
vol. 1, p. 8　'reduce 2130/4004 hundredths[21]

約　　之　　乃　　二千二　　之　　一千六十五
*yuē*　*zhī*　*nǎi*　*èr qiān èr*　*zhī*　*yī qiān liù shí wǔ*
reduce　3OBJ　then　[2][$10^3$][2]　*zhī*　[1][$10^3$][6][10][5]
to get 1065/2002 [hundredths]'[22]
(this result is subsequently expressed as 0.0053196)

Moreover, Li Zhizao used the terms *líng* and *yòu* to concatenate the integer and fractional parts of a mixed number; this configuration can be seen twice with *líng* in (3) and once with *yòu* in (4) where the tenths rank *fēn* is treated in the same way as a measure word. Both linking terms can also be used to concatenate numbers of different units or decimal ranks; example (3) is with *yòu* and (5) uses *líng*.

(3)　十斤　　零　　五分　斤　之　二
in *Tóng wén*　*shí jīn*　*líng*　*wǔ fēn*　*jīn*　*zhī*　*èr*
*suàn zhǐ*　[10] *jīn*　*líng*　[5] *fēn*　*jīn*　*zhī*　[2]
vol. 1, p. 9　'10 2/5 *jīn*

又　　七兩　　零　　二分　之　一
*yòu*　*qī liǎng*　*líng*　*èr fēn*　*zhī*　*yī*
*yòu*　[7] *liǎng*　*líng*　[2] *fēn*　*zhī*　[1]
and 7 1/2 *liǎng*'

---

[21] The word *lí* 鏊 stands for 0.01. Decimal ranks smaller than 1 were expressed using a scale starting with the words *fēn* 分 for $10^{-1}$, *lí* 鏊 for $10^{-2}$, *háo* 毫 for $10^{-3}$, *sī* 絲 for $10^{-4}$, *hū* 忽 for $10^{-5}$, *wēi* 微 for $10^{-6}$; for each digit of the decimal part, the rank was specified. For example, 0.123 would be expressed as [1] *fēn* [2] *lí* [3] *háo*. The linguistic pattern "Num + Rank name" was the same as for units of measurement.

[22] The expression of the first fraction follows "Denominator's name + *fēn* + MW + *zhī* + Numerator's name", for the second it is "Denominator's name + *zhī* + Numerator's name".

(4)  　　　　七錢　　八分　　又　　七分　　分　　之　　二
in *Tóng wén*　*qī qián*　*bā fēn*　*yòu*　*qī fēn*　*fēn*　*zhī*　*èr*
 *suàn zhǐ*　[7] *qián*　[8] tenth　*yòu*　[7] *fēn*　tenth　*zhī*　[2]
vol. 1, p. 8　'7 *qián* [and] 8 tenths [of a *qián*] and 2/7 tenth'

(5)  　　　　一千丈　　　零　　四分　　三釐
in *Tóng wén*　*yī qiān zhàng*　*líng*　*sì fēn*　*sān lí*
 *suàn zhǐ*　[1][10³] *zhàng*　*líng*　[4] tenth　[3] hundredth
vol. 1, p. 8　'1000 *zhàng* and 4/10 3/100' (i.e. 1000.43 *zhàng*,
　　　　　expressed a few line further as 100,043 hundredths)

In Contemporary Chinese, the linking term *líng* can still be used outside the linguistic numeration to connect compounds expressed in a scale of currency units as in (6) and (7), and in the year-month scale in (8) and (9) where no zero is involved[23].

In (6) we have the succession of ranks of the currency unit scale *yuán* (unit of currency), *jiǎo* (tenth of a *yuán*) and *fēn* (hundredth[24] of a *yuán*) which are all needed to express the price 3.85 in which there is no void rank, and actually the final *fēn* can be omitted. But to say 3.05 as in (7), the word *jiǎo* is not used and the linking term *líng* comes between the compounds "Num + *yuán*" and "Num + *fēn*", and *fēn* can again be omitted without causing any ambiguity.

(6)  三　　元　　八　　角　　　　　　五　　分
　　*sān*　*yuán*　*bā*　*jiǎo*　　　　*wǔ*　*fēn*
　　[3]　yuan　[8]　tenth of a yuan　[5]　hundredth of a yuan
　　'3 yuan and 85 cents'

(7)  三　　元　　零　　五　　分
　　*sān*　*yuán*　*líng*　*wǔ*　*fēn*
　　[3]　yuan　*líng*　[5]　hundredth of a yuan
　　'3 yuan and 5 cents'

In (8) and (9) we can see the time-unit scale formed by the two words *nián* [year] and *yuè* [month]. There is no other unit expected between them, but the idiomatic pattern still requires the linking term *líng*.

---

[23] Incidentally, such a similarity of treatment between the series of pivots and scales of measure words should be noted.

[24] Historically, *fēn* was a designation of tenths within a scale of decimal values each 1/10 of the preceding; the terms of this series are nowadays used as the prefixes for the International Systems of Units. But the same words were also used in a scale decreasing by a factor 1/100; the first term *fēn* designated hundredths; this meaning is still extant in currency units.

(8) 一 年 零 五 個 月
    *yī   nián   líng   wǔ   ge   yuè*
    [1]   year   *líng*   [5]   CLF   month
    'one year and five months'

(9) 一 年 零 十一 個 月
    *yī   nián   líng   shíyī   ge   yuè*
    [1]   year   *líng*   [10][1]   CLF   month
    'one year and eleven months'

The same word *líng* is also used to designate zero in Contemporary Chinese, both the number and the symbol to write it, but nothing suggests that the *líng* used within the morphosyntax of number names is a linguistic transposition of the zeros used in symbolic notation (and indeed not all zeros of a number written in Arabic digits have a *líng* counterpart in the number name); it is a term used to link a [*digit*][*pivot*] compound to another similar compound or to the units digit.

In Contemporary Chinese no pivot name automatically calls for the item *líng*; it occurs when there is a *gap* within the series of pivots in the number name, in other words when the next pivot to be said is not the next smaller one in the list of pivots. That is to say *líng* is obligatory: *i)* after the inner pivots *qiān* [$10^3$] and *bǎi* [$10^2$] if the following digit is not associated with *bǎi* [$10^2$] or *shí* [10] respectively; *ii)* after the outer pivots *yì* [$10^8$] and *wàn* [$10^4$] if the following digit is not associated with *qiān* [$10^3$].

The five configurations which trigger the use of *líng* with numbers of three or four digits are visualized in Table 2. For example the three-digit number 105 is expressed as [1][$10^2$] *líng* [5]. Two zeros in a row in the positional notation of a four-digit number correspond to only one *líng*, for example 1001 is expressed as [1][$10^3$] *líng* [1]; in other words there is no one-to-one mapping of the number name components onto the digits in positional notation.

Tab. 2: Configurations of numbers of three or four digits which trigger the use of *líng*

| Configurations | | | | | Name of number in |
|---|---|---|---|---|---|
| *thousand* | *hundred* | *ten* | *unit* | | Contemporary Chinese |
| ≠0 | =0 | ≠0 | ≠0 | → | [*thousand*][$10^3$]*líng*[*ten*][10][*unit*] |
| ≠0 | =0 | ≠0 | =0 | → | [*thousand*][$10^3$]*líng*[*ten*][10] |
| ≠0 | ≠0 | =0 | ≠0 | → | [*thousand*][$10^3$][*hundred*][$10^2$]*líng*[*unit*] |
| ≠0 | =0 | =0 | ≠0 | → | [*thousand*][$10^3$]*líng*[*unit*] |
|  | ≠0 | =0 | ≠0 | → | [*hundred*][$10^2$]*líng*[*unit*] |

To illustrate what happens at the level of the outer pivots, let us consider the number 1,305,000,080 but let us see it as 13,0500,0080, cut into slices of four digits consistent with the use of the two outer pivots *wàn* [$10^4$] and *yì* [$10^8$]. These will not be immediately followed by a [*digit*][$10^3$] block because the next four-digit slices start with zeros shaded here: 13,0500,0080. This number is expressed as shown in (10). Again there is no one-to-one mapping of the occurrences of *líng* onto the zeros in positional notation: the two series of zeros shaded above each trigger one occurrence of *líng*, and the three zeros left unmarked do not.

(10) 十 三 億 零 五 百 萬 零 八 十
*shí sān yì líng wǔ bǎi wàn líng bā shí*
[1] [3] [$10^8$] *líng* [5] [$10^2$] [$10^4$] *líng* [8] [10]
'1,305,000,080'

This is the situation in *correct* Contemporary Chinese, but Zhejiang speakers[25] can drop *líng* after an outer pivot (never after an inner pivot) even though they know it is considered *incorrect*.

In order to reach the present-day use of *líng* in names for integers, there is no doubt that some standardization was implemented at some point; this may have occurred during the late 19th or early 20th centuries. In the mathematics books of the *Bái fú táng suàn xué cóng shū* collection[26] edited from 1872 to 1877 during the Qing dynasty, some authors like Shi Yuechun used the same number names as Qin Jiushao in the 13th century, possibly as a conscious revival of the Song dynasty mathematical tradition[27], whereas some other authors do not use *líng* in number names. The shift from free-coined phrases and unplanned linguistic creation to *approved* standardized expressions required going through a process that Haugen (1983) called *corpus planning*. Haugen distinguished four steps: *selection of norm* (which is societal and exterior to the language); the *codification of the norm*; *implementation of function* (includes the activities of writers and institutions); and the *elaboration of function* (involves the production of a linguistic corpus complying with the norm). To complete this part of the history of Contemporary Chinese integer names would require more research on the definition and implementation of standard number names in the late 19th and/or early 20th centuries.

---

[25] Information provided in Paris separately by some thirty speakers who are occasional speakers of the Wenzhou dialect but have Mandarin (standard Contemporary Chinese) as their major communication language at home and as their first educational language; they also declared not to know how to express large numbers in Wenzhou dialect.

[26] The *Bái fú táng suàn xué cóng shū* 白芙堂算學叢書 collection is composed of 23 books of mathematics edited by Ding Quzhong 丁取忠 (1810-1877); more details can be found in Wu Wenjun (2000, 200-203).

[27] The 14th century saw developments in the calculations with an abacus and lost interest in some domains of mathematics explored earlier and related to calculations with counting rods. Chinese mathematicians regained interest in Qin Jiushao's work only after the introduction of European mathematics in the 16th and 17th centuries.

## 2-2. Changes concerning the use of [1]

The script of Shang and Zhou inscriptions concatenates the transcription of a sequence [*digit*][*pivot*] forming only one character, thus making it impossible to know whether the morpheme [1] was used before pivots.

The *Suàn shù shū* (early 2$^{nd}$ century BCE) contains hundreds of integers written in the Chinese language, allowing thorough comparisons of all possible configurations. The morpheme [1] is used before all pivots of a number name but the highest one, with no exception. Please check shaded *yī* 一 [1] in (11), (12), (14) and (16) for examples of the former situation, and shaded in (13), (15) and (16) for the latter position.

(11)            二百      一十  
in *Suàn shù shū*    *èr bǎi*    *yī shí*  
strip 172        [2][10$^2$]   [1][10]  
                 '210'

(12)            二千      一十      六  
in *Suàn shù shū*    *èr qiān*   *yī shí*    *liù*  
strip 20         [2][10$^3$]   [1][10]    [6]  
                 '2016'

(13)            錢        ∅百       五十  
in *Suàn shù shū*    *qián*    ∅*bǎi*     *wǔshí*  
strip 76         *qián*    ∅[10$^2$]    [5][10]  
                 '150 *qián*' (*qián* 錢 is a currency unit)

(14)            七千      一百      二十      九  
in *Suàn shù shū*    *qī qiān*   *yī bǎi*   *èr shí*   *jiǔ*  
strip 176       [7][10$^3$]   [1][10$^2$]   [2][10]   [9]  
                 '7129'

(15)            ∅千      八十      九  
in *Suàn shù shū*    ∅*qiān*   *bā shí*   *jiǔ*  
strip 172       ∅[10$^3$]   [8][10]   [9]  
                 '1089'

(16)            ∅萬      一千      五百      二十      銖  
in *Suàn shù shū*    ∅*wàn*    *yī qiān*   *wǔ bǎi*   *èr shí*   *zhū*  
strip 47        ∅[10$^4$]   [1][10$^3$]   [5][10$^2$]   [2][10]   *zhū*  
                 '11520 *zhū*' (*zhū* 銖 is a unit of weight)

The *Jūyán Xīnjiǎn* 居延新简 bamboo strips were excavated in 1974 and date from the 1st century CE[28]. The use of [1] with pivots in this corpus, as in (17), is identical to what can be seen in the *Suàn shù shū*.

| (17) | | 凡 | ∅萬 | 一千 | 一百 |
|---|---|---|---|---|---|
| in *Jūyán Xīnjiǎn* | | *fán* | ∅*wàn* | *yī qiān* | *yī bǎi* |
| 4454: E.P.T53:129 | | total | ∅[$10^4$] | [1][$10^3$] | [1][$10^2$] |
| | | 'a total of 11100' | | | |

There were three texts of mathematics excavated at Dunhuang. They are *Pelliot chinois 2667*[29], *Pelliot chinois 3349*[30] and *Stein 19 Recto*[31]. The date range of their composition spans from the 1st to the 10th century CE, which cannot be narrowed down further. These texts exhibit some changes concerning the use of the morpheme [1]: it is used even with the highest pivot, as in [1][$10^2$] in (18), except if this pivot is [10] in which case [1] is optional; compare (18) without [1] to (19) with [1].

| (18) | 二五 | 如 | ∅十 | 自相乘得 | 一百 |
|---|---|---|---|---|---|
| in *Pelliot chinois* | *èr wǔ* | *rú* | ∅*shí* | *zìxiāng chéng dé* | *yī bǎi* |
| *3349* | [2][5] | as | ∅[10] | REF multiply get | [1][$10^2$] |
| | '2 times 5 is 10, which multiplied by itself gives 100' | | | | |

| (19) | 一十一 | 萬 | 五千 |
|---|---|---|---|
| in *Pelliot chinois* | *yī shí yī* | *wàn* | *wǔ qiān* |
| *3349* | [1][10][1] | [$10^4$] | [5][$10^3$] |
| | '115000' | | |

Actually in *Pelliot chinois 3349*, [1] is absent not only before [10] in [10][$10^4$] but also before [$10^2$] and [$10^3$] in the expressions [$10^2$][$10^4$] and [$10^3$][$10^4$]. These three number names follow the pattern [*digit*][*pivot*] with the pivot [$10^4$] and the digit-slot being occupied by [10], [$10^2$] and [$10^3$] respectively, revealing a different behavior when these inner pivots are multiplicands of the outer pivot [$10^4$] than when they are used as pivots.

The *Nine Chapters* [*Jiǔ Zhāng Suàn Shù* 九章算術] is a text originally written during the Han dynasty (206 BC–220 AD), but the known edition was prepared in the 7th century CE[32] and might have undergone linguistic amendments. In this transmitted corpus, the morpheme [1] is used before the highest pivot in a number name, even if it is [10], even before [10][10⁴] the compound of an inner and an outer pivot [33]; this hardly argues for a generalization of this feature before the 7th century CE (and the anteriority of Dunhuang texts of mathematics) because some features of the *Nine Chapters* might simply reflect the choices of its authors and not the ordinary linguistic situation at the time of writing.

---

[28] Accessed in January 2012 on the website of the Academia Sinica (Taiwan); there are 5 relevant instances.

[29] Accessed in May 2012 in Gallica (Bibliothèque nationale de France); there are 4 relevant instances.

[30] Accessed in May 2012 in Gallica (Bibliothèque nationale de France); there are 19 relevant instances.

[31] Accessed in May 2012 on the website of the International Dunhuang Project; there are 3 relevant instances.

[32] Chemla and Guo Shuchun (2004: 43-46).

The 13th century mathematician Qin Jiushao followed the norm set by the editors of the *Nine Chapters*. The situation is still the same in Contemporary Chinese, but when the highest pivot of a number name is [10] it is not obligatorily preceded by [1][34]; this double capacity makes [10] a boundary point between digits and pivots.

## *2-3. Two morphemes for 2*

Zhou Shengya (1984) explains that *liǎng* in Old Chinese was not used like other cardinal numbers; it was used for objects naturally coming in pairs (e.g. *liǎng ěr* 兩耳 'both ears', *liǎng shǒu* 兩手 'both hands') or in historical names like *Liǎng Zhōu* 兩周 'Western and Eastern Zhou'; this was somewhat akin to English *both*. Before a pivot *liǎng* showed some similarity to the English "*a couple of + Noun*" in its approximate meaning *some* or *a few*. Only *èr* could be used in exact number names. But in Contemporary Chinese, the two numerical morphemes *èr* and *liǎng* are encountered in names for integers.

However, *liǎng* can never replace *èr* as an ordinal number: Contemporary Chinese can use almost all cardinal names as ordinals with or without the prefix *dì* 第; the only exception is *liǎng*. With or without the prefix *dì*, only *èr* can be used to state the second ordinal position as in (20) for "*the second floor*", whereas the cardinal form of 2 before a classifier is usually *liǎng* as in (21) for "*two floors*".

(20)   他      住       二       層
       *tā*    *zhù*    *èr*     *céng*
       3SG    dwell    [2]      floor
       'He lives on the second floor.'

---

[33] One instance only in the main text of the section 5-10 (not in the commentary).

[34] In Contemporary Chinese, the ordinary way of expressing 10 is *shí* [10]. The expression *yī shí* [1][10] occurs when extra clarity is intended as for example when stating accounts or voicing calculations. In any configuration, before *shí* [10] the usual shift of *yī* [1] to falling tone before a rising tone syllable is neutralized and the pronunciation of [1] remains *yī* with a high tone.

14(21)  他    住     一套    兩       層     的     房子
     *tā*   *zhù*   *yī tào*   *liǎng*   *céng*   *de*   *fángzi*
     3SG   dwell   [1] CLF   [2$_{variant}$]   floor   PART   apartment
     'He lives in a two-floor apartment.'

Moreover *liǎng* is preferred over *èr* as a cardinal before a classifier. Only the classifier *liǎng* (50 g) favors *èr* for reasons of euphony; please compare (22) and (23). With other classifiers *èr* is possible, just less common and more formal.

(22)  兩      個     人
     *liǎng*   *ge*   *rén*
     two   CLF   person
     'two persons'

(23)  二     兩       餃子
     *èr*   *liǎng*   *jiǎozi*
     two   CLF   dumpling
     'two *liǎng* of dumplings' (i.e. 100 g of dumpling)

In Contemporary Chinese the two morphemes *liǎng* and *èr* can occur in exact number names before the pivots [$10^2$], [$10^3$], [$10^4$] and [$10^8$]. But only *èr* is used with [10], and in the unit-slot when there are other digits above. A search I made in June 2012 on the search engine Baidu provides the distribution given in Tab. 3 of the two morphemes before each pivot. The item *liǎng* is more frequent than *èr* and the frequency increases with higher pivots. Higher pivots are treated as classifiers; only [10] is not, and as above for its behavior with respect to the digit [1], it holds a special position among the series of pivots.

Tab. 3: Distribution of *èr* 二 and *liǎng* 兩 with pivots on *Baidu* 百度 (June 2012)

| Pivots | *shí* [10] | | *bǎi* [$10^2$] | | *qiān* [$10^3$] | | *wàn* [$10^4$] | | *yì* [$10^8$] | |
|---|---|---|---|---|---|---|---|---|---|---|
| Compounds | *èr shí* | *liǎng shí* | *èr bǎi* | *liǎng bǎi* | *èr qiān* | *liǎng qiān* | *èr wàn* | *liǎng wàn* | *èr yì* | *liǎng yì* |
| Number of occurrences | 100,000,000 | 440,000 | 42,700,000 | 48,500,000 | 25,500,000 | 71,400,000 | 14,000,000 | 56,300,000 | 1,680,000 | 11,200,000 |
| Occurrences with *liǎng* in percentage of the total with *èr* and *liǎng* | 0.4 % | | 53 % | | 74 % | | 80 % | | 87 % | |

In a CCL search (accessed in January 2012) for *liǎng* in complex number names (isolated *liǎng* [*pivot*] sequences were rejected since they could be approximate numbers as explained above), the earliest occurrence[35] is found in a text first published in 1343, and there are four instances dated to the turn of the 19th and 20th centuries[36]. Then the situation changes dramatically with 6712 occurrences in the contemporary corpus (after 1911). It follows that the use of *liǎng* in exact number names is a phenomenon which started in the late 19th or early 20th centuries and developed further in the 20th century. The earlier occurrences look like anachronisms, yet they could be isolated instances of a potential linguistic novelty.

In any case, Contemporary Chinese can now freely choose between the morphemes *èr* [2] and *liǎng* [2$_{variant}$] to express the digit 2 with all pivots but [10]; for example the number 2222 can be said as in (24) or (25):

(24)  二千　　二百　　二十二
　　　*èr qiān*　*èr bǎi*　*èr shí èr*
　　　[2][10$^3$]　[2][10$^2$]　[2][10][2]
　　　'2222'

(25)  兩千　　兩百　　二十二
　　　*liǎng qiān*　*liǎng bǎi*　*èr shí èr*
　　　[2$_{variant}$][10$^3$]　[2$_{variant}$][10$^2$]　[2][10][2]
　　　'2222'

Native speakers claim that it is more *correct* to use only *èr* for all occurrences of the digit 2 when reading a number, i.e. outside of any syntactic or contextual incorporation. This assumed *correctness* surely does not imply that they favor the variants with *èr* exclusively.

---

[35] In *Sòng shǐ – Zhì dì yī bǎi sān shí wǔ – Shí huò xià sì* 宋史•志第一百三十五•食貨下四, a history of the Song dynasty first published in 1343 under the Yuan dynasty. The transmitted text gives *mǐ jià dàn liǎng qiān wǔ bǎi zhì sān qiān* 米價石兩千五百至三千, that is 'the price for husked grain is 2500 to 3000 per *shí* [a unit of capacity]' with 2500 expressed as [2$_{variant}$][10$^3$][5][10$^2$]. The same text contains 53 instances of *èr qiān wǔ bǎi* 二千五百, i.e. [2][10$^3$][5][10$^2$] with *èr*. The text in CCL is likely a 1934 edition published by the Shànghǎi Shāngwù Yìnshūguǎn 上海商務印書館, and the occurrence of *liǎng* here might be an editorial error.

[36] One instance in *Guī lú tán wǎng lù* 歸廬譚往錄 by Xu Zongliang 徐宗亮 (1828-1904), one in *Kāngxī xiáyì zhuàn* 康熙俠義傳 and two in *Xù Jì Gōng zhuàn* 續濟公傳 both published under Emperor Guangxu 光緒 (1875-1908).

# 3. THE MORPHOSYNTAX OF LINGUISTIC NUMERATIONS AND THE SYNTAX OF QUANTIFICATION

If we identify the compounds formed from one digit and one pivot in Chinese and other languages as a quantification noun phrase, then we assume the digit and multiplicative pivot to have different semantic and syntactic functions: the latter is construed as a noun, a measure word or a classifier, whereas the former works as a quantifier. Some features of number names in Chinese confirm this approach and show that the morphosyntax of number names is rooted in the syntax of quantification. But other characteristics can diverge from what such an analysis should imply, and they reveal that the structure of number names possesses some degree of autonomy with regard to the patterns used to express quantification.

## *3-1. Order of the elements expressing* digit *x* pivot

The relative order of the two elements in the compounds expressing the products *digit* x *pivot* is the same as the order of the quantification pattern in Chinese. This is also true in English, but it is not true in Tibetan numeration which also relies on [*digit*][*pivot*] compounds, whereas the order of the quantification pattern is "*Noun + Num*"[37]. The order [*multiplicative pivot*][*digit*] could have been possible; it exists in Iraqw (Tanzania), Ndom (Papua New Guinea) and Yorùbá (Nigeria). Yet Tibetan ordinary quantification order with [*digit*] after [*pivot*] can resurface when expressing round numbers whether by direct juxtaposition with *khir* [$10^4$] and *bum* [$10^5$] or by inserting *phrag* after *bgya* [$10^2$] and *stong* [$10^3$] [38]. The incorporation of the resulting expressions into more complex number names is limited and requires a conjunction.

## *3-2. Conjunctions between compounds expressing* digit *x* pivot

In languages like Chinese and English which use the [*digit*][*pivot*] order, the occurrence of a reverse sequence [*pivot*][*digit*] indicates that two compounds [*digit*][*pivot*] have been concatenated by direct juxtaposition and that the sum of those two compounds is implied; there is no risk of confusion and no pragmatic need for a conjunction. A conjunction was nevertheless used in Chinese before the

---

[37] For more information about the structure of quantification phrases in various Tibeto-Burman languages, one can refer to Xu Dan (2010) and the chapter by Fu Jingqi in the present book.
[38] Goldstein et al. (1991: 199). Wylie's transliteration is used for Tibetan words.

3rd century BCE, for example in [10] *yòu* [5] for 15. The use of linking terms to concatenate compounds expressing the products *digit* x *pivot* is a feature which shows that the morphosyntax of a linguistic numeration is rooted in the syntax of noun phrases, or more precisely in the patterns used to concatenate noun phrases. Striking examples are found in the oracular inscriptions of the Shang dynasty (14th to 11th century BCE): the phrase in (26) is an expression of "*fifteen dogs*" which does not exhibit an unbreakable number name for "*fifteen*", but rather a succession of two quantification phrases "*ten dogs*" and "*five dogs*" linked by *yòu*.

(26)　　　　十　　犬　　有(又)　　五　　犬
in H32775　　*shí*　*quǎn*　*yòu*　　*wǔ*　*quǎn*
　　　　　　[10]　dog　　and　　　[5]　dog
　　　　　　'fifteen dogs'

The expressions "[10] Noun + *yòu* + [5] Noun" and "[10] Noun + *yòu* + [5]" with nouns other than *dog* are also encountered in Shang inscriptions. Eventually the disappearance of the noun after [10] made it possible for the compound [10] *yòu* [5] to occur, and produced a number name independent from its context of syntactic incorporation.

Finally even if the use or disuse of a conjunction began because of the implementation of some language planning, they are accepted and have been transmitted because they fit into the syntax of noun phrases.

### *3-3. Similarity of Chinese pivots with classifiers and nouns*

Contemporary Chinese possesses the two morphemes *èr* and *liǎng* to express the digit 2. The latter is favored before classifiers and also before all multiplicative pivots but [10]. This means that all pivots but [10] in the sequence [*digit*][*pivot*] bear a syntactic resemblance to classifiers.

In Contemporary Chinese *an investment of five hundred million* can be expressed with the ordinary numeral $[5][10^8]$ as in (27), or using a "Num + CLF + Noun" phrase where the numeral is the digit [5] and the noun is the pivot $[10^8]$ as in (28).

(27)　　五億　　　　　的　　　投資
　　　*wǔ yì*　　　　*de*　　*tóuzī*
　　　$[5][10^8]$　　PART　investment
　　　'an investment of five hundred million'

(28)　　五個億　　　　的　　　投資
　　　*wǔ ge yì*　　　*de*　　*tóuzī*
　　　[5] CLF $[10^8]$　PART　investment
　　　'an investment of five hundred million'

The latter expression is a round number which cannot be incorporated in a more complex number name; it nevertheless reveals that the highest pivot [$10^8$] can be readily reinterpreted as a noun.

## *3-4. Use of [1] with multiplicative pivots*

The compound [1][*pivot*] parallels the quantification pattern with [1] as the quantifier and [*pivot*] as the quantified item. This is consistent with the idea of an isomorphism between a [*digit*][*pivot*] sequence and a quantification phrase.
Contemporary Chinese follows the same pattern [1][$10^3$][1][$10^2$] as English *one thousand one hundred* to express 1100, with [1] before each pivot. But in the *Suàn shù shū* (beginning of the 2[nd] century BC) it was [$10^3$][1][$10^2$], and in French it is *mille cent* [$10^3$][$10^2$], both number names illustrating the fact that a linguistic numeration can deviate from the quantification pattern.
The use of [1][*pivot*] is a sign of such an isomorphism. Beware however that the absence of [1] with a pivot is not necessarily a proof of independence between the two sub-systems. For example, Arabic uses morphological means to mark plurality; the word for $10^3$ has the three forms *alf*, (singular), *alfān* (dual) and *ālāf* (plural). The morpheme [1] cannot appear before the singular *alf*; the isomorphism between number names and quantification phrases is nevertheless established by the morphology of plurality[39].

## *3-5. Elliptic number names*

The relative autonomy of number name systems allows the production of *elliptic number names* which are not meant to be incorporated in quantification phrases. For example, in Contemporary Chinese, the name for 150 can be the regular *yī bǎi wǔ shí*, i.e. [1][$10^2$][5][10], or the elliptic *yī bǎi wǔ*, i.e. [1][$10^2$][5] dropping the last pivot [10]. The latter number name is not ambiguous because 105 is pronounced *yī bǎi líng wǔ*, i.e. [1][$10^2$] *líng* [5]. However before a classifier, and also before the outer pivots [$10^4$] or [$10^8$], the elliptic form is rejected as awkward, even though it is unambiguous, and only the complete form with the last pivot can be incorporated. Elliptic forms which drop a pivot are known in other languages. For example, in French, to express a price of 2,500,000 one can use the abbreviated *deux millions cinq*, i.e. [2][$10^6$][5], instead of *deux millions cinq cent mille*, i.e. [2][$10^6$][5][$10^2$][$10^3$], and nobody, in the context of stating a price, would mistake it for the number 2,000,005 which is also pronounced *deux millions cinq*, i.e. [2][$10^6$][5][40]. This elliptic French expression cannot be followed by the designation of a currency; its ambiguity is perhaps not the sole obstacle to its incorporation since this kind of restriction is observed even with unambiguous Contemporary Chinese expressions.

---

[39] But inconsistencies in the marking of plurality within number names are in turn a sign of the relative autonomy of their morphosyntax with respect to the syntax of quantification. For example, English pivots used in exact number names never bear the plural marker, while they do when used as approximate numbers: compare *three thousand* [3][$10^3$] (an exact number name) and *thousands* [$10^3$]$_{PL}$ (a round number expressed with *thousand* analyzed as a countable noun).

The utterance of a pivot provides an informative context leading the listener to infer that the digit which immediately follows should be the number of units of the rank just below the one which was said previously, for example [5] in Contemporary Chinese [1][$10^2$][5] is expected to belong to the rank 10 which is just below $10^2$. In the same manner [5] in French [2][$10^6$][5] is expected to be the number of units at the rank $10^5$, unless another rank is specified. In other words, in the process of extracting numerical information step by step from the sequence of morphemes composing the number name, a listener (resp. speaker) will anticipate what should follow[41], thus allowing elliptic numerical expressions in which the last rank name is dropped; however, this is possible if and only if the digit pronounced previously belongs to the rank just below the previously uttered rank. But when putting together a quantification expression and preparing to use a classifier, noun or measure word after the numerical expression, the last digit should express the number of units of the noun or measure word. If this is not what is intended, the name of the last rank cannot be omitted, and Chinese has [1][$10^2$][5][10] before a classifier regardless of the fact that [1][$10^2$][5] is an unambiguous expression of 150 in Contemporary Chinese.

## 4. CONCLUSION

The Chinese numeration system was always decimal-based, with the same digits [1] through [9], and a series of multiplicative pivots including [10], [$10^2$], [$10^3$], [$10^4$] in all ancient sources and with the addition of [$10^8$] which is the highest commonly agreed pivot in Contemporary Chinese. Contemporary Chinese numeration does not possess special terms for teens, which are expressed as [10][*unit digit*] or [1][10][*unit digit*], nor for tens, which are expressed as [*tens digit*][10].

Complex Chinese number names are concatenations of compounds made of a *digit* followed by a *pivot*; these compounds basically follow the quantification patterns "Num + Noun", "Num + MW" and "Num + CLF". Historically, these forms employed with simple numerals (say 1 to 10 in Chinese) were existing patterns readily used by analogy for the development of a more complex linguistic numeration, pivots being equivalent to nouns, measure words or classifiers.

---

[40] This example was proposed by Robert Iljic (CRLAO, France).
[41] This phrasing is somehow naïve, I do not assume to understand the underlying cognitive processes.

Nevertheless the autonomy of the morphosyntax which produces number names with regard to the syntax of quantification manifests itself in Contemporary Chinese with elliptic number names. These drop the last pivot, and cannot be incorporated before classifiers, where they would conflict with the encoding and decoding processes corresponding to the linguistic pattern of quantification.

Historical changes occurred in the use of the linking terms *yòu* and *líng* and in the how the digit [1] was used with pivots. The expressions of 105 and 150 at various periods exhibited in Tab. 4 illustrate the effects of these changes. This historical description provides criteria to help determine the date of composition for excavated texts or of re-editions for transmitted texts. Taking into account the rules followed by the number name system at a given time can also help reconstruct damaged fragments.

Tab. 4: Evolutions of the names for 105 and 150 in Chinese

|  | Chinese names for 105 | Chinese names for 150 |
|---|---|---|
| 13th-3rd centuries BCE: the conjunction *yòu* is attested but was used irregularly, whether [1] was used with [$10^2$] is inaccessible. | [$10^2$] *yòu* [5] or [1][$10^2$] *yòu* [5], [$10^2$][5] or [1][$10^2$][5]. | [$10^2$] *yòu* [5][10] or [1][$10^2$] *yòu* [5][10], [$10^2$][5][10] or [1][$10^2$][5][10]. |
| At the beginning of the 2nd century BCE in the *Suàn shù shū*. | [$10^2$][5]. | [$10^2$][5][10]. |
| After a change between the 1st and the 7th centuries concerning the use of [1]. | [1][$10^2$][5]. | [1][$10^2$][5][10]. |
| With the introduction of the linking term *líng* during the 12th and 13th centuries; later disused, eventually revived in the late 19th or early 20th centuries. | [1][$10^2$] *líng* [5]. | [1][$10^2$][5][10], and the elliptic [1][$10^2$][5] dropping the last pivot is also possible in Contemporary Chinese with no ambiguity with the name for 105 which requires *líng*, however it cannot be incorporated in a quantification phrase. |

The generalization of the use of [1] before the first pivot to appear in a number name is a change more likely to have been provoked by language planning than because of the intrinsic evolution of the language. This feature had reached its full development in the 7th century CE edition of the *Nine Chapters*, a 7th century edition of a Han dynasty text. Nevertheless there is still much uncertainty concerning the whole situation from the 1st to 10th centuries CE.

Contemporary Chinese names for integers over 100 use the linking term *líng* in a manner irreconcilable with any one-to-one mapping of the components of the number names onto positional notation but allowing dropping of the last pivot with no ambiguity. This was already the case in Qiu Jiushao's 13th century book, but lost in Li Zhizao's 17th century writings; he used *líng* to link various quantities but not within names for integers. We do not know who initiated the revival of 13th century number names in the late 19th or early 20th centuries. This normative action of language planning was able to succeed because its use of the linking tem *líng* was consistent with the ordinary concatenation of noun phrases. However, some speakers still tend to drop *líng* after outer pivots.

Finally, since the 20th century, some numbers have developed free variants due to the possibility of choosing between the morphemes *èr* and *liǎng* before any pivot besides 10; this is a grass-roots development running against the regularity imposed by language planning.

The features of Chinese linguistic numeration and their evolution reveal the tension between a tendency to shorten number names versus an inclination to maintain their morphosyntax within the framework of the syntax of quantification and of noun phrases, thus minimizing the variety of cognitive processes involved in encoding and decoding. The coexistence of these opposing mechanisms is likely to prove valid cross-linguistically. Extrapolating from the case of Chinese, one can suggest two characteristics of linguistic numerations which are highly susceptible to historical changes: *i*) whether one times a pivot is merely [*pivot*] or requires saying [1]; *ii*) whether compounds meaning a product *digit* x *pivot* are directly juxtaposed or linked with conjunctions. In addition, one of these two characteristics can undergo changes without the other one being affected; and as these features are ruled by conflicting evolutionary tendencies and language planning, changes are not necessarily mono-directional or irreversible.

## Corpora and Dictionaries